\documentclass{amsart}
\usepackage{color}

\newtheorem{thm}{Theorem}

\theoremstyle{definition}
\newtheorem{definition}[thm]{Definition}
\newtheorem{example}{Example}

\newtheorem{prop}{Proposition}[section]
\newtheorem{remark}[thm]{Remark}

\newtheorem{problem}{Problem}
\newtheorem{theorem}{Theorem}
\newtheorem{lemma}{Lemma}

\usepackage{amsmath}

\newcommand{\blem}{\begin{lemma}}
\newcommand{\elem}{\end{lemma}}
\newcommand{\bexer}{\begin{exe}}
\newcommand{\eexer}{\end{exe}}
\newcommand{\beq}{\begin{eqnarray}}
\newcommand{\eeq}{\end{eqnarray}}
\newcommand{\bthm}{\begin{theorem}}
\newcommand{\ethm}{\end{theorem}}
\newcommand{\bex}{\begin{example}}
\newcommand{\eex}{\end{example}}

\newcommand{\bdefi}{\begin{definition}}

\newcommand{\edefi}{\end{definition}}

\newcommand{\bprop}{\begin{prop}}
\newcommand{\eprop}{\end{prop}}

\newcommand{\bpf}{\begin{proof}}
\newcommand{\epf}{\end{proof}}

\def\be{\begin{equation}}
\def\ee{\end{equation}}

\newcommand{\U}{\mathbb{U}}

\numberwithin{equation}{section}

\input amssym.def
\input amssym

\begin{document}
\author{M. Mateljevi\' c}
\address{Faculty of mathematics, University of Belgrade, Studentski Trg 16,
Belgrade, Serbia} \email{miodrag@matf.bg.ac.rs}
\title[Schwarz lemma]{Schwarz lemma, and  Distortion   for harmonic functions via length  and area}

\date{28 Dec,2016}
\maketitle

This is very rough working  version (Version  3,   4/27/2018).
\section{Introduction and Basic definitions}
We give sharp estimates  for distortion of  harmonic  by means  of area and length of the corresponding surface.
In 2016 \cite{RgSchw1}(a),   the author has posted  the  current Research project  Schwarz lemma, the Carath\'{e}odory and Kobayashi Metrics and Applications in Complex Analysis.\footnote{motivated by  S. G. Krantz  paper \cite{krantz}}  Various discussions regarding the subject  can also  be found  in the Q\&A section on Researchgate under the question  '' What are the most recent versions of The Schwarz Lemma ?'',\cite{RgSchw1}(b).
During the fall semester 2017 at Belgrade seminar \cite{Bg_sem},   we  have  communicated about  Schwarz lemma  and we have posted the
arXiv  paper \cite{MMSchw_Kob}, in which  we  have
considered  various version  of Schwarz lemma and its relatives related to harmonic and holomorphic  functions including  distortion of  harmonic mappings,
and  several variables.  For  the results of  \cite{MMSchw_Kob}  see also  \cite{MMSchw_Kob1}.
For example,  in Section  \ref{ss_schw} we prove several  optimal   versions  of planar Schwarz  lemma for   real   valued  harmonic maps $h$
from $\mathbb{U}$ into $I_0=(-1,1)$ ( see Theorem \ref{burg1} and \ref{mmsch1}, related to  the case  $h(0)=a$, $a\in I_0$;  and   Theorem \ref{mmsch11}  and  \ref{th_khav1}, for the case   $f(a)=b$, $a\in \mathbb{U}$).
In particular  if $a=0$
a part of Theorem   \ref{burg1} is reduced  to  classical Schwarz  lemma for harmonic maps.\\
Note that  Theorem  \ref{th_khav1}  yields  solution of  D. Khavinson   extremal problem for harmonic functions in  planar case, cf. \cite{Khav,Maz}.\\
From  Theorem   \ref{burg1} we also  derive Theorem   \ref{burg2} which is a version  of planar Schwarz  lemma  for  complex   valued  harmonic maps $h$
from $\mathbb{U}$ into  itself,   and a  version of   the boundary   Schwarz  lemma,  see  Theorem   \ref{zhu1}.

During my work on the subject, D. Kalaj  gave an interesting communication,   cf. \cite{Kal_Sim17}(from which we have learned about his  arXiv  papers  \cite{Kal_jga,Kal_minimal}), and immediately we have realized that we can adapt our previous consideration to connect with  his work.
In particular,  using different    approach  we
can give new insight to   these results as well as further results.

We first need some definitions.
\bdefi \label{def1}a1) By $\mathbb{C}$ we denote   the complex plane,  by   $\mathbb{U}$   the   unit disk  and by  $\mathbb{T}$ the   unit circle. \\
a2) In planar case  $G\subset \mathbb{C}$  with euclidean norm the notation $\Lambda_f(z)$, $z\in G$,   is used instead of  $||(df)_z||$.
For a function $h$, we use notation  $D_1h=h_{x}^{\prime}$ and   $D_2h=h_{y}^{\prime}$  for partial derivatives;
$\partial h=\frac{1}{2}(h_{x}^{\prime}-ih_{y}^{\prime })$ and $\overline{\partial}h=\frac{1}{2}(h_{x}^{\prime }+ih_{y}^{\prime })$; we also use
notations $Dh$ and $\overline{D}h$ instead of $\partial h$ and
$\overline{\partial }h$ respectively when it seems convenient.
We use the notation  $\lambda_f (z)= \big||\partial f (z)|-|\bar\partial f (z)|\big|$ \, and   $\Lambda_f (z)= |\partial f (z)|   +
|\bar\partial f (z)|,$   if $\partial f(z)$   and  $\bar\partial
f(z)$ exist.\\
a3)   If   $f$ is  harmonic   in   $\mathbb{U}$   then we write   $f= g+\overline{h}$, where $g$ and $h$ are analytic
in   $\mathbb{U}$. If  $|g'(z)|\geq |h'(z)|$, then   $\Lambda_f (z)=|g'(z)| + |h'(z)|$,  $\lambda_f(z)=|g'(z)| - |h'(z)|$  and therefore
$\Lambda_f (z)  + \lambda_f(z)= 2 |g'(z)|$. In general,      $\Lambda_f (z)  + \lambda_f(z)=2 \max\{|g'(z)|,|h'(z)|\}$.
\edefi
\bdefi  \label{def2}
b1)  For   a  $C^1$ mapping   $u: \mathbb{U} \rightarrow \mathbb{R}^m$,
set  $S= u(\mathbb{U})$, $D[u]=\int_{\mathbb{U}} (|D_1u|^2 + |D_2u|^2 )dx dy$,  $E=|D_1 u|^2$,  $G=|D_2 u|^2$, $F=D_1 u \cdot D_2 u$,     $J_u=\sqrt{EG -F^2}$, and     $A=A(S)=A(u)=\int_{\mathbb{U}}J_u dx dy$.\\
b2) We say that  $u$ is $\underline{K}$-qc if   $E+F \leq 2\underline{K} J_u$; in planar case  this definition  is a small  modification  of  the standard definition of coefficient of quasi-conformality, see  Remark \ref{rmQc3}.

b3)  Suppose  that      $u: \mathbb{U} \rightarrow \mathbb{R}^m$  is   harmonic  on $\mathbb{U}$.
Then  $ u={\rm Re} F$, where $F$ is analytic.   Set $\underline{D}[F]=\int_{\mathbb{U}} |F'(z)|^2 dx dy$.

b4) If $f$ is a function on $\mathbb{T}$,  we associate to $f$ a curve  $\gamma=\gamma_f$ defined by
$\gamma(t)= f(e^{it})$, $t\in[0,2\pi]$, and  we denote   by   $L=L(f)=  |\gamma_f|$  the  length of   $\gamma_f$.\\
\edefi
It is easy to check that \\
$D_1 u ={\rm Re} F'$  and   $D_2 u =-{\rm Im} F'$  and \\
(A)  $|D_1 u|^2 + |D_2 u|^2 = |F'|^2 $, and therefore\\
$D[u]=\underline{D}[F]$.\\
If in addition  $u$ is conformal at some point,   $|D_1 u|= |D_2 u|$  and therefore\\
$|{\rm Re} F'|= |{\rm Im} F'|$  and    (ii)  $|F'|^2= 2 |D_1 u|^2$.

We will use the following hypothesis in the  sequel \\
{\rm ($H_m$):} $u: \mathbb{U} \rightarrow \mathbb{R}^m$   is  harmonic  on $\mathbb{U}$,  and    $S=u(\mathbb{U})$,\\
{\rm($\underline{H}_m$)}:  In addition to  {\rm ($H_m$)}  we suppose that \\
{\rm  (h1):}  $u$  is   continuous on  $\overline{\mathbb{U}}$  and   $\gamma=\gamma_u$  is rectifiable.\\
In  Section    \ref{se_area}  we prove:\\
{\rm (I0)}  If  $u$ satisfies  \rm($\underline{H}_m$),
then  {\rm (i1)}:  $4 \pi A \leq L^2$, where   $A=A(u)$  and    $L=L(u)$. \\
{\rm(I1)}  If in addition to   {\rm ($H_m$)}   we suppose that $u$ satisfies\\
{\rm(h2):}  $A(S)$  is finite  and \\
{\rm(h3):}  $u$ is  $\underline{K}$-quasiconformal,\\
then
$$D[u]=\underline{ D}[F]=\pi \bigl(\sum_{k=1}^\infty (k |\hat{F}(k)|^2) \bigr) \leq 2 \underline{K }A(S).$$
{\rm I2)}  If in addition to   {\rm ($H_m$)}   and   {\rm(h2)}  we suppose that   {\rm(h4):} $u$ is  conformal, then
$$ A=\int_{\mathbb{U}} |D_1u|^2 dx dy = \frac{\pi}{2} \bigl(\sum_{k=1}^\infty (k |\hat{F}(k)|^2) \bigr).$$
{\rm d2)}  In particular,    $\pi \Lambda^2_u(0)\leq D[u]$  with equality iff (ii) $\gamma=\gamma_u$  is a circle  given by  $u_k=a_kx -b_ky$, $k=1,2,3,...,m$,   where   $|a|= |b|$  and $a\cdot b=0$.\\
Hence using the isoperimetric inequality,\\
{\rm d3)}   $2\pi \Lambda_u(0)\leq L$   with equality iff (ii).\\
In \cite{Kal_jga} a version of  (d3)  is proved.

In  Section    \ref{Se_dia} we first consider the cases when $m=2,3$.

For  $f: \mathbb{U} \rightarrow \mathbb{R}^2$ we use notation  $G=f(\mathbb{U})$.
For convenience of the reader we also  first   suppose that $u$ is harmonic on   $\overline{\mathbb{U}}$. In this case $L$  is length of  $\partial G$.
In  Section    \ref{Se_dia} using  Proposition \ref{propMM0}, we consider  distortion  of  harmonic    functions on $\mathbb{U}$  related to  diameter $dia(G)$  of image domain $G=f(\mathbb{U})$, see Theorem \ref{length0}.   Then use an  inequality between     $dia(G)$  and
the length $L(G)$   of boundary of $G$  we estimate  distortion  via  $L(G)$ and  prove
$2\pi  \lambda_f(0) \leq L$,  see Theorem \ref{length1}.
In  \cite{Kal_Sim17}   a version of  the part (b) of  Theorem \ref{length1} is proved for diffeomorphisms.\\
Under  the hypothesis \rm($\underline{H}_m$),  $m\geq 2$,
it is convenient  to introduce for given  $z_0\in \mathbb{U}$   the tangent plane $Z=Z_{\mathbf{y}^0}$, $\mathbf{y}^0= u(z_0)$. Then we use the projection $p$   onto  $Z$
and apply planar result on   $p\circ u$  to prove   $2\pi (1 - |z_0|^2) \lambda_u(z_0) \leq L$, see Theorem \ref{length2}  and Theorem \ref{tH0}.  If in addition $f$ is  conformal at $z_0$, then  $\Lambda_u(z_0)=\lambda_u(z_0)$
and the previous inequality holds with  $\Lambda_u(z_0)$ instead of  $\lambda_u(z_0)$, see Theorem \ref{length3}.\\

In  Section  \ref{ss_gen} we outline a proof of  Theorem \ref{t_incr}. Using this result one can  show that some  of  the above described results  hold under more general   hypothesis then {\rm($\underline{H}_m$)}(see for example Theorem \ref{length+0}).\\

These results are communicated in November 2017,\cite{Bg_sem}.
\section{Schwarz  lemma}\label{ss_schw}
For a hyperbolic plane domain $D$, we   denote by $\rho_D$(or  $\lambda_D$) the hyperbolic density  and by abusing notation the  hyperbolic metric occasionally.
\blem
If  $G$ and $D$ are simply connected domains  different from $\mathbb{C}$    and \\
$\omega \in {\rm  Hol}(G,D)$, then
$\rho_D(\omega z)|\omega'(z)|\leq \rho_G(z) $,  $z\in G$  and
$$\rho_D(\omega z,\omega z')\leq \rho_G(z,z'), \quad z,z'\in G .$$
\elem
We denote the right half plane by  $\Pi$.

\bprop
If   $\omega$   is holomorphic from  $\Pi$  into itself, then
$$|\omega'(z)|\leq   \frac{{\rm Re} \omega(z)}{{\rm Re} z}.$$
If in addition  $\omega$  maps $\mathbb{R}^+$   into itself, then   $ |\omega'(1)|\leq {\rm Re} \omega(1)= \omega(1)$  and therefore   $\omega'(1)\leq \omega(1)$.
\eprop
\bdefi
By $\mathbb{C}$ we denote   the complex plane  by   $\mathbb{U}$   the   unit disk  and by  $\mathbb{T}$ the   unit circle.
For   $z_1 \in \mathbb{U}$,   define  $$T_{z_1}(z)= \frac{z-z_1}{1-\overline{z_1} z},$$    $\varphi_{z_1}
= -T_{z_1}$.

d1)  Throughout   this paper  by   $\mathbb{S}(a,b)$ we denote the set  $(a,b)\times \mathbb{R}$,   $-\infty \leq  a < b \leq \infty$,  and in particular we write   $\mathbb{S}_0$ for $\mathbb{S}(-1,1)$.
Note that   $\mathbb{S}(a,b)$ is a strip if   $-\infty < a < b < \infty$  and   $\mathbb{S}(a,+\infty)$ is a half-plane if  $a$ is a real number,  and   $\mathbb{S}(-\infty,+\infty)=\mathbb{C}$.  By  $\lambda_0$  and  $\rho_0$  we denote hyperbolic metrics  on  $\mathbb{U}$  and  $\mathbb{S}_0$  respectively.

d2)  Set  $I_0=(-1,1),\quad  \mbox{ and for}\quad   a\in I_0 \quad  \mbox{define}$
$$s=s(a)=\tan(\frac{\pi}{4}(a +1)),\quad     e=e(a)=\cot\big(\frac{\pi}{4}(a +1)\big),\quad    \mbox{and}$$
 $$X(r)=X^+(r,a)=\frac{4}{\pi} \arctan(s \frac{1+ r}{1-r})-1,\quad X^-(r,a)=1- \frac{4}{\pi} \arctan\left(e \frac{1+ r}{1-r} \right), z\in \mathbb{U}.$$
Further it is convenient to introduce the functions $A,B,A_s$ and $B_s$    by      $A(r)= (1+r) (1-r)^{-1},  B(r)= (1-r) (1+r)^{-1} $, $A_s (r)=s A(r), B_s(r)=s B(r)$,  and     $Y(r)= X^+(|z|,|a|)$.

d3)  Set   $c=(a+1)/2$,    $\overline{c}= 2\pi c$,   $\alpha=\alpha(c)= \alpha(a)=\overline{c}/2=(a+1)\pi/2$.

It is convenient to write  $f_y(x)=f(x,y)$.

(A0)   It is straightforward  to check
$$X_a^-(r)= \frac{4}{\pi} \arctan(B_s(r))-1,\quad   X_a^-(r) \leq a \leq X_a^+(r),$$
$X^+(r,a)$ (respectively  $X^-(r,a)$)     is increasing (respectively  decreasing) in both  variables  $r$  and $a$,   $X_1^+ =1$  and $X_{-1}^- = -1$.

\edefi
Note that  $s=s(|a|)=\tan(\frac{\pi}{4}(|a| +1)) \quad    \mbox{for}  \quad   a\in  \mathbb{U}$.

Since    $X(r)= \frac{4}{\pi} (\arctan\circ A ) - 1$  and   $A'_s(r)= 2s (1-r)^{-2}$, we find
\beq \label{burgD}
X'(r)=\frac{4}{\pi} \frac{2 s (1-r)^{-2}}{1+ A_s^2(r)}=    \frac{4}{\pi}   \frac{2s}{(1-r)^2 +s^2 (1+r)^2},   0\leq r < 1.
\eeq
In a similar way  since  $X_-(r)=   \frac{4}{\pi} (\arctan\circ B)  -1$  and   $B'_s(r)= - 2s (1+r)^{-2}$, we find
\beq \label{burgDb}
 X_-'(r)= -\frac{4}{\pi}   \frac{2s}{(1+r)^2 +s^2 (1-r)^2},     0\leq r < 1.
\eeq

Next
$$X(0)= \frac{4}{\pi} \arctan\left(\tan \frac{\alpha(c)}{2}\right) -1  =\frac{4}{\pi} \frac{\alpha(c)}{2}-1  =\frac{4}{\pi} (a+1)\frac{\pi}{4} -1   =a,$$
and  by (\ref{burgD}),
\beq \label{burgD0}
X'(0)=\frac{4}{\pi} \frac{2s }{1+s^2}=\frac{4}{\pi} \frac{2\tan \frac{\alpha(c)}{2} }{1+  \left(\tan \frac{\alpha(c)}{2}\right)^2}=\frac{4}{\pi} 2  \tan (\alpha/2) \cos^2(\alpha/2)  =\frac{4}{\pi} \sin \alpha,
\eeq
and in a similar way  using   (\ref{burgDb})
\beq \label{burgD0b}
\quad X_-'(0)= - X'(0) =-\frac{4}{\pi} \sin \alpha.
\eeq

Suppose that $f$ is harmonic map from $U$ into  $I_0=(-1,1)$   with $h(0)=a$.
Using    a version of Schwarz lemma  \cite{MMSchw_Kob1}, we will show
\begin{equation}
\rho_0(fz,a)= |\ln \frac{s(fz)}{s(a)}|  \leq  \ln\frac{1+r}{1-r}, z\in \mathbb{U} .
\end{equation}
This inequality  is equivalent to  $X^-(|z|,a) \leq  f(z)\leq  X(|z|)=X^+(|z|,a), z\in \mathbb{U}$.

\bthm\label{burg1}  If $u_1,u_2\in (-1,1)$, then
\begin{equation}\label{eq:dist0}
\rho_0(u_1,u_2)=|\ln \frac{s(u_2)}{s(u_1)}|.
\end{equation}
Let  $h$  be  a  real   valued  harmonic map
from $\mathbb{U}$ into $I_0=(-1,1)$   with $h(0)=a$, $a\in I_0$. Then
\beq \label{burgEst1}
X^-(|z|,a) \leq  h(z)\leq  X(|z|)=X^+(|z|,a), z\in \mathbb{U},  \\
\mbox{and} \quad  |(dh)_0| \leq  X'(0)=\frac{4}{\pi} \sin \alpha.\label{eq:har1}
\eeq
\ethm
If $a=0$, then $a_1= \tan \frac{\pi}{4} =1$  and   $X(|z|,0) = \frac{4}{\pi} \arctan |z|$. Hence we get classical Schwarz  lemma for harmonic maps
which states  $|h(z)|\leq  X(|z|)=X(|z|,0)=\frac{4}{\pi}\arctan |z|$.

\bpf
We use
$$\sec \theta =\frac{1}{\cos \theta},\quad  \tanh^{-1}z= \frac{1}{2}\ln\frac{1+z}{1-z},\quad
I(x)= \int_0^x \sec t\, dt =2\tanh^{-1}(\tan \frac{x}{2})$$  and
\begin{equation}\label{eq:dhyp-strip0a}
    \rho_{0}(w)= {\rm  Hyp}_{\mathbb{S}_0}(w)= \frac{\pi}{2}  \frac{1}{\cos (\frac{\pi}{2} u)},\mbox{for}\quad  w\in \mathbb{S}_0, \quad \mbox{where}\quad u={\rm Re} w.
\end{equation}

Since $A(\tan \frac{x}{2})= \tan (\frac{\pi}{4} + \frac{x}{2})$,  we find  $I(x)=\ln \big(\tan (\frac{\pi}{4} + \frac{x}{2})\big)$.

If $u_1,u_2\in (-1,1)$, $u_1\leq u_2$, then  using the change of variables  $t=\frac{\pi}{2} u$,  $t_k=\frac{\pi}{2} u_k$, $k=1,2$,
we have
\begin{equation}
\rho_0(u_1,u_2)= \frac{\pi}{2} \int_{u_1}^{u_2} \frac{du}{\cos (\frac{\pi}{2} u)}= \int_{t_1}^{t_2} \frac{du}{\cos (t)} =I(t_2)-I(t_1),
\end{equation}
and therefore  since   $I(t_k)=\ln  \big(\tan (\frac{\pi}{4} + \frac{x}{2})\big)= \ln s(u_k)$, $k=1,2$, we find
\begin{equation}
\rho_0(u_1,u_2)= \ln\frac{\tan\frac{\pi}{4}(u_2+1)}{\tan\frac{\pi}{4}(u_1+1)}.
\end{equation}
Hence   (\ref{eq:dist0})  follows.

If $h(0)=a$  and recall  we set  $s=s(a)=\tan(\frac{\pi}{4}(a +1))$,    by  a version of Schwarz lemma  \cite{MMSchw_Kob1}, for  $z\in \mathbb{U}$ we find

\beq \label{burgEst1a}
\tan(\frac{\pi}{4}(h(z) +1)) \leq s \frac{1+|z|}{1-|z|},\quad \mbox{ie.}\quad   h(z)\leq  X(|z|)=X(|z|,a).
\eeq
XX   If $a=0$, then $s(0)= \tan \frac{\pi}{4} =1$  and   $X(|z|,0) = \frac{4}{\pi} \arctan |z|$.
If we set
$g=-h$,    then  $g(0)= -a$,  and by  (\ref{burgEst1a}), we find    $- h(z)\leq  X(|z|)=X(|z|,-a)$, ie.  $ h(z)\geq - X(|z|,-a)$.  Hence  one can derive  (\ref{burgEst1}).

But,  we prefer the following approach.
If $z\in K_r$, then
\begin{equation}
|\ln \frac{s(fz)}{s(a)}|  \leq  \ln\frac{1+r}{1-r} .
\end{equation}
We can rewrite this inequality as
\begin{equation}
\frac{s(fz)}{s(a)} \leq  \frac{1+r}{1-r}\quad \mbox{if}\quad  s(a)\leq s(fz), \quad \mbox{and}\quad  \frac{s(a)}{s(fz)} \leq  \frac{1+r}{1-r} \quad \mbox{if}\quad  s(a)\geq s(fz) .
\end{equation}

Hence  if  $f(z)\geq a $, we find  $s(fz)\leq  A_s(r)$  and therefore  $f(z) \leq  X(r)$.
Further it is convenient to introduce       $B(r)=s (1-r) (1+r)^{-1} $, and     $X_-(r)=   \frac{4}{\pi} (\arctan\circ B)  -1$.
Hence  if  $fz\leq a $, we find  $s(a)\leq s(fz)  A_s(r)$  and therefore  $f(z) \geq  X_-(r)$.

Hence, by  (A0),    we find  $ X_-(r) \leq  f(z)\leq  X(r)$.

Next by (\ref{burgD})  and  (\ref{burgDb}),  we have
$$X'(r)= \frac{4}{\pi}   \frac{2s}{(1-r)^2 +s^2 (1+r)^2},\quad  X_-'(r)= -\frac{4}{\pi}   \frac{2s}{(1+r)^2 +s^2 (1-r)^2},     0\leq r < 1,$$
$\mbox{and in particular  by  (\ref{burgD0})  and (\ref{burgD0b})},    X'(0)=\frac{4}{\pi} \sin \alpha,\quad X_-'(0)=-\frac{4}{\pi} \sin \alpha$,
and therefore
(\ref{eq:har1})  follows.
\epf
XX After writing  the previous version we have  realized that the inequality   (\ref{burgEst1}) in Theorem  \ref{burg1}  is covered by   \cite{bu}, but our proof is completely  different.

\bdefi
d1)  For $a\in (-1,1)$, let ${\rm Har}^a$    denote the family of all real  valued  harmonics maps
$f$ from $\mathbb{U}$ into $(-1,1)$ with $f(0)=a$.

d2)   For $a\in \mathbb{U}$  and $b\in (-1,1)$, set   $L(a,b)= L(a,b) = \sup|(du)_a|$, where the supremum is taken over all  real  valued  harmonics maps $u$
from $\mathbb{U}$ into $(-1,1)$  with $u(a)=b$.

d3)  For $a\in \mathbb{U}$  and  $\ell \in  T_a\mathbb{C}$  a unit vector,
set  $L(a)= \sup|(du)_a|$  and $L(a,\ell)= \sup|(du)_a(\ell)|$, where the supremum is taken over all  real  valued  harmonics maps
from $\mathbb{U}$ into $(-1,1)$.
\edefi

Now,  we can restate and strength the part of  Theorem  \ref{burg1}:
\bthm\label{mmsch1}
If  $a\in (-1,1)$  and   $h\in {\rm Har}^a$, then
\begin{equation}\label{burgEst2}
(i) \quad h(z)\leq  X(|z|),   \quad  (ii) \quad   |(dh)_0| \leq  X'(0)=\frac{4}{\pi} \sin \alpha  \quad   \mbox{and}\quad (iii) \quad  L(0,a)=\frac{4}{\pi} \sin \alpha(a).
\end{equation}
\ethm
\bpf
We need only to prove (iii). There is  a      conformal mapping $f$   of  $\mathbb{U}$  onto  $\mathbb{S}_0$  with  $f(0)=a$   and  $f'(0)>0$; then for  harmonic function    $u_0= {\rm Re} f$  the equality  holds in (iii).
\epf

\bthm\label{mmsch11}
Let  $h$  be  a  real  valued  harmonics map  from $\mathbb{U}$ into $(-1,1)$ with $f(a)=b$, $a\in \mathbb{U}$. Then
 \be
h(z) \leq \frac{4}{\pi} \arctan\left(\frac{1 + |\varphi_a(z)|}{1 - |\varphi_a(z)|} \tan \frac{\alpha(|b|)}{2}\right) -1,
 \ee
\be\label{mmsch2}
|(dh)_a| \leq  \frac{4}{\pi}\frac{\sin \alpha(|b|)}{1- |a|^2} .
\ee
\ethm
\bpf Set  $w= \varphi_a(z)$.   Apply Theorem \ref{mmsch1}  on  $h^a=h\circ \varphi_a$, we find   $h^a(z)\leq X(|z|)$.  Hence  $h(w)= h^a(z)\leq  X(|\varphi_a(w)|)$.
Since we can identify   $(d\varphi_a)_0$  with  $1-|a|^2$,    using   $(dh^a)_0= (dh)_a\circ (d\varphi_a)_0$   and Theorem \ref{mmsch1}  we prove    (\ref{mmsch2}).
\epf

Further set
$$A_0(z)=\frac{1+z}{1-z},\quad   \mbox{and   let}\quad    \displaystyle \phi = i \frac{2}{\pi} {\mathrm {ln}} A_0;$$
that is
$\phi= \phi_0\circ A_0$,  where   $\displaystyle{\phi_0 = i \frac{2}{\pi} {\mathrm {ln}}}$.
Let $\hat{\phi}$ be defined by   $\displaystyle{\hat{\phi} (z)=-\phi (iz)}$. Note that $\phi$  maps $I_0=(-1,1)$ onto y-axis  and
$\hat{\phi}$ maps $I_0$   onto itself.

If  $\hat{u}= {\rm Re } \hat{\phi}$, then
\begin{equation}\label{eext1}
    \hat{u}=\frac{2}{\pi} {\rm arg} \left(\frac{1+iz}{1-iz}\right)
\end{equation}
and   $\hat{u}$  maps $I_0=(-1,1)$ onto itself.

Let  $a\in (0,1)$  and  $\ell \in  T_a\mathbb{ C}$. There is  a    conformal mapping $f=f_\ell$  of  $\mathbb{U}$  onto  $\mathbb{S}_0$  with  $f(a)=0$   and  $f'(a)\ell>0$.
We will show that     $u=u_\ell= {\rm Re} f_\ell$  is extremal.
In particular,   there is  a      conformal mapping $f$   of  $\mathbb{U}$  onto  $\mathbb{S}_0$  with  $f(a)=0$   and  $f'(a)>0$; set   $u_0= {\rm Re} f$.
\bthm \label{th_khav1}
If  $a\in (-1,1)$   and  $\ell \in  T_a\mathbb{C}$, then
\begin{enumerate}
\item $L(a)= (u_0)'_r (a)=\frac{4}{\pi}(1- |a|^2)^{-1}$  and
\item  $L(a,\ell)= L(a)=(du_\ell)_a(\ell)= \frac{4}{\pi}(1- |a|^2)^{-1}$.
\end{enumerate}
\ethm
This yields  solution of  D. Khavinson   extremal problem for harmonic functions in  planar case, cf. \cite{Khav,Maz}.
\bpf   (1)  By  hypothesis   $\rho_0 (f(a)) |f'(a)|=  2 (1- |a|^2)^{-1}$,  $\rho_0 (f(a))=\rho_0(0)= \frac{\pi}{2}$,  $(u_0)'_r (a)=f'(a)$  and $|(du_0)_a|=|\nabla u_0 (0)|=
\frac{4}{\pi}(1- |a|^2)^{-1}$.

(2)  Recall  there is  a    conformal mapping $f=f_\ell$  of  $\mathbb{U}$  onto  $\mathbb{S}_0$  with  $f(a)=0$   and  $f'(a)\ell>0$.  If    $u=u_\ell= {\rm Re} f_\ell$, then
$(du)_a(\ell)= {\rm Re} \big(f'(a)\ell  \big)$.
We leave the interested  reader to fill details.
\epf
\bthm \label{th_khav2}
Let  $h$  be  a  real  valued  harmonics map
 from $\mathbb{U}$ into $(-1,1)$ with $f(a)=b$, $a\in \mathbb{U}$. Then
\be\label{mmsch3}
(i)\quad |(dh)_a| \leq   \frac{4}{\pi}\frac{\sin \alpha(|b|)}{1- |a|^2}, \quad  (ii)\quad   L(a,b)=\frac{4}{\pi}\frac{\sin \alpha(|b|)}{1- |a|^2} .
\ee
\ethm
\bpf
There is  a    conformal mapping   of  $\mathbb{U}$  onto  $\mathbb{S}_0$  with  $f(a)=b$. We leave the interested  reader to  show that   $u_0= {\rm Re} f$
is extremal  for (i)  and therefore  (ii) holds.
\epf
\subsection{Schwarz lemma at the boundary}
\bthm\label{burg2}
Let  $h$  be  a  complex    valued  harmonic map
from $\mathbb{U}$ into itself    with $h(0)=a$, $a\in \mathbb{U}$. Then
$$
 |h(z)|\leq  X(|z|)=X^+(|z|,|a|), z\in \mathbb{U}.
$$
\ethm
\bpf
Using   rotation around $0$   and  Theorem \ref{burg1} one can prove this result.
\epf

\begin{thm}\label{zhu1}
 Let   $f: \mathbb{U} \rightarrow \mathbb{U}$    be  harmonic  and $s=s(f(0)$.   Further assume that there is a point  $b \in \mathbb{T}$ so that $f$ extends continuously to $b$, $|f(b)| = 1$ {\rm(}say that $f(b) = b'${\rm)}, and
$f$ is   $\mathbb{R}$- differentiable at $b$.
Then
$$|\Lambda_f(b)| \geq \frac{2}{s \pi}.$$
\end{thm}
\bpf By   (\ref{burgD}), we find

$$\lim_{r\rightarrow 1_-} X'(r)= \frac{2}{s \pi }.$$
The rest of proof is based on  Theorem  \ref{burg2}  and the following proposition.
\epf
We leave the interested reader to  prove the following propositions:
\begin{prop}\label{pr_g0}
{\rm (a)} Let   $f: \mathbb{U} \rightarrow \mathbb{U}$.
 Assume that  there is a point  $b \in \mathbb{T}$ so that $f$ extends continuously to $b$, $|f(b)| = 1$ {\rm(}say that $f(b) = c${\rm)}, and
and $f$ is  $\mathbb{R}$- differentiable at $b$.

{\rm (b)} Further  assume that  there is a function $A$  such that  $A:[0,1] \rightarrow [0,1]$, $A'(1)$ exists and  $M_f(r) \leq A(r)$.

Then  $|\Lambda_f(b)|\geq |f'_r (b)| \geq A'(1)$.
\end{prop}
\bpf
Without loss  of  generality   we can suppose  that  $c=b=1$. By (b),
$$|\frac{f(1)- f(r)}{1-r}|\geq \frac{1- M_f(r)}{1-r} \geq B(r):= \frac{1- A(r)}{1-r}.$$
Hence if $r\rightarrow 1_- $, we have   $|f'_r (b)| \geq A'(1)$. Since by definition of   $\Lambda_f(b)$,  $\Lambda_f(b)\geq |f'_r (b)|$ it  completes proof.
\epf
\begin{prop}
Under the above hypothesis, if  there exists $f'(b)$, then

{\rm(i)}  $|f'(b)| \geq A'(1)$.
\end{prop}
\section{Distortion  of  harmonic  functions related to  diametar  and   length}\label{Se_dia}
We advise the reader to recall Definition \ref{def1}.
In \cite{kavu} and  \cite{MMSchw_Kob}  in particular, it is proved (see also \cite{abr}, Theorem 6.16, Proposition  6.19, cf.  \cite{Khav,bu,Maz}):
\begin{prop} \label{propMM0}
If  $u$ is a  harmonic map from $\U$  into $I_0=(-1,1)$, then
\beq\label{ineqMM0}
|\nabla u (0)| \leq  \frac{4}{\pi}.
\eeq
\end{prop}
For convenience of the reader we outline a proof.
Throughout   this paper  by   $\mathbb{S}(a,b)$ we denote the set  $(a,b)\times \mathbb{R}$,   $-\infty \leq  a < b \leq \infty$,  and in particular by  $\mathbb{S}_0=\mathbb{S}(-1,1)$.
The mapping  $f_0$  defined by    $f_0(w)= \tan(\frac{\pi}{4} w)$       maps $\mathbb{S}_0$ onto  ${\mathbb U}$.

If we denote by  $\rho_0$   hyperbolic density on  $\mathbb{S}_0$,  then  using $f_0$  we can check  that  for $w=u+iv \in \mathbb{S}_0$,
\begin{equation}\label{eq:dhyp-strip0a}
    \rho_{0}(w)= {\rm  Hyp}_{\mathbb{S}_0}(w)= \frac{\pi}{2}  \frac{1}{\cos (\frac{\pi}{2} u)}.
\end{equation}
It is known from the standard course of Complex Analysis  that  there is an   analytic function  $\omega$ on $\U$   such that    $u = {\rm Re}\omega$ on  $\U$.
Since   $\omega$ is holomorphic map from  $\U$ into   $\mathbb{S}_0$, then by a very special case  of  Schwarz-Ahlfors-Pick  lemma(see also  the property (I)),

\begin{equation}\label{eq:Ah-Sch1a0}
    \rho_0 (\omega(z)) |\omega'(z)|\leq  2 (1- |z|^2)^{-1},\quad z\in \U,
\end{equation}
where $\rho_0$  is given by  (\ref{eq:dhyp-strip0a}).

Since  $\frac{\pi}{2}\leq \rho_{0}(w)$  and  $|\omega'|=|\overline{\nabla u}|=  |\nabla u|$,  we have   $(\ref{ineqMM0})$.

In particular,  if  $\omega$  is a    holomorphic function   from the   unit disk  $\mathbb{U}$ into $\mathbb{S}_0$  with $\omega(0)=0$,
we have
$|\omega'(0)|\leq \frac{4}{\pi}$  with the  equality   iff $\omega$ is  a conformal mapping  of  $\mathbb{U}$ onto  $\mathbb{S}_0$.
\begin{remark}
Note that one can derive  Theorem \ref{mmsch1}  from   (\ref{eq:Ah-Sch1a0}). Namely,  by the above notation   $\rho_0 (u(z)) |\nabla u(z)|\leq  2$.
\end{remark}
\bdefi
c1)  If $g$ is a holomorphic function on  $\mathbb{U}$  by $\hat{g}_k$ we denote its Taylor coefficient and write
$g(z)= \sum_{k=0}^\infty \hat{g}_k z^k $. Note that  $k! \hat{g}_k=  g^{(k)}(0)$.
c2)  For a set $M \subset \mathbb{R}^m$   by   ${\rm d=dia}(M)$ we denote  the  diameter of $M$.

\edefi
\begin{thm}\label{length0}
Let $f= g+\overline{h}$ be  complex valued harmonic   in   $\mathbb{U}$  which satisfies  {\rm(H1)}. \\
Then \\
{\rm(i)} $\pi \Lambda_f(0)\leq 2d$.\\
{\rm(ii)} $2d \leq  L$.\\
For function $u_d(z)= \frac{d}{\pi }\arg \frac{1+z}{1-z}$  the  equality holds in {\rm(i)}.
\end{thm}
It is interesting that  $p(z)=d \cdot x/2$  is not extremal for the inequality   (i).
\bpf
We can suppose that  $f(0)=0$   and using rotations   that   $\Lambda_f(0) =|\hat{e}_1|$ and   $\hat{e}_1= df_0(e_1)=k e_1$, where  $k=\Lambda_f(0)$.
Set  $p(w)= u$  and $F= p \circ f$.  Then  $p(G)$  is an interval of length equal or  less then $d$,
and  by Proposition  \ref{propMM0},   $\pi |\nabla F(0)|\leq 2d$. Since  $|\nabla F(0)|=\Lambda_f(0)$,  we get the first inequality  of  (i).
We   leave to the reader    to show that  $2d \leq  L$.\\
\epf
\bdefi
\begin{itemize}
\item[d1)]  For $a\in\mathbb{ R}$, define $S_a=\{w: {\rm Re}w < a\}$   and   let    $\mathcal{P}_a$  denote the family of all functions  $f$  holomorphic   in   $\mathbb{U}$  for which
$f(\mathbb{U}) \subset S_a$.
\item[d2)] If $H$ is a holomorphic function  on $\mathbb{U}$ which   has  zero at $0$ at least of order $2$  and $a  \in \mathbb{C}$,  it is straightforward to check that  there are unique holomorphic functions $g=g_H$ and $h=h_H$ on $\mathbb{U}$  such that

\item[{\rm (h0)}:]  XX   $g(0)=h(0)=0$  and   $g'= -  H + a$, $h'= z^{-2}H$.
\end{itemize}
Note that  in this setting $a=g'(0)$  and set  $f_{H,a}= g_H +i \overline{h_H}$.  If it is not confusing we write $f_H$ instead of  $f_{H,a}$ and    in this way we associate a unique  harmonic function $f_H$   to $H$.
We  say that $H$ satisfies
\begin{itemize}
\item[($h^\prime_1$):]  if it   satisfies  {\rm (h0)}  with  $2 H \in \mathcal{P}_1$, \\
 and that
$f=g +\overline{h}$  satisfies
\item[(\underline{h}0)]   with respect to    $H$:  if  $H$  satisfies    ($h^\prime_1$).
\item[d3)]  It is convenient to say that $f=g +\overline{h}$  satisfies
\item[ (\underline{h}1)]  with respect to  $\nu$  if:
  $g(0)=h(0)=0$  and
\beq\label{def_Belt0}
g'=  \frac{1}{1+ z^2\nu} \quad   {\rm and} \quad h'= \frac{\nu}{1+ z^2\nu},
\eeq
where
\item[(i3):]  $\nu=\omega z^{-2}$  and  $\omega\in {\rm Hol}(\mathbb{U},\mathbb{U})$   has  zero at $0$ at least of order $2$.

\item If    $\nu$ satisfies  (i3) there a unique  $f$, which we denote by  $f^{\underline{\nu}}=f_0^\nu$,   such that
$g'$ and   $h'$  are given  by  (\ref{def_Belt0}).
\end{itemize}
\edefi
Note if    $g'$  is  given  by  (\ref{def_Belt0}), then  $g'(0)=1$,
and if   $\nu$ satisfies  (i3), then by an application of classical Schwarz lemma,  $|\nu(z)|\leq 1$,  $z\in \mathbb{U}$,  and the function  $(1+ \omega)^{-1}$   (defined by  $z\mapsto (1+ z^2\nu)^{-1}$)  is holomorphic on
$\mathbb{U}$.

We leave to the interested reader to show   that   $f=g +\overline{h}$  satisfies (\underline{h}1)  with respect to  $\nu$   iff  it  satisfies (\underline{h}0)  with respect to   $H$ with  $H=\frac{\nu z^2}{1+ \nu z^2}$.

\begin{thm}\label{length1}
Let $f= g+\overline{h}$ be complex valued  continuous on  $\overline{\mathbb{U}}$  and   harmonic  on $\mathbb{U}$  which satisfies  {\rm($\underline{H}_2$)}.
\footnote{$\gamma(t)= f(e^{it})$, $t\in[0,2\pi]$
is a rectifiable curve and  $L= |\gamma|$ is length of   $\gamma$.}
Then
\begin{itemize}
\item [{\rm a)}]  $2 \pi k |\hat{g}_k|\leq L$, $k\geq 0$. In particular
\item [{\rm a1)}]  $2 \pi |g'(0)|\leq L$,   with equality in the case  $a= g'(0)>0$  iff
\item [{\rm ($h_a$)}:]  $g'= -  H + a$, $h'= z^{-2}  H$,  and   $2 H \in \mathcal{P}_a$,  where  $H$ has  zero at $0$ at least of order $2$.
\item [{\rm a2)}]    $2 \pi \max\{|g'(0)|,|h'(0)|\}\leq L$
\item [{\rm a3)}]    $2 \pi |g'(0)|\leq L$    with equality iff
\item [{\rm(i4)}:]   $f=cf^\nu +c_1$, where  where    $\nu$  satisfy  {\rm(i3)}.
\item [{\rm b)}]   $2 \pi (1-|z|^2)|g'(z)|\leq L$, $z\in \overline{\mathbb{U}}$  with equality iff     {\rm(i5)}:  $f=c f^{\underline{\nu}} \circ  \varphi_z +c_1$,
where $c,c_1\in \mathbb{C}$  and     $\nu$  satisfy  {\rm(i3)}.
\item [{\rm b1)}]  $2 \pi (1-|z|^2)|\lambda_f(z)|\leq L$, $z\in \overline{\mathbb{U}}$.
\end{itemize}
\end{thm}
As a corollary   we get,   $\pi (|g'(0)|+ |h'(0)|)\leq L$ and  since  $\Lambda_f (z)  + \lambda_f(z)=2 \max\{|g'(0)|,|h'(0)|\}$,
$\pi (|\Lambda_f (0)  + \lambda_f(0)|)\leq L$  and    $2\pi  \lambda_f(0) \leq L$.

Set $iX(t)=f'_t e^{-it}$.  If $g'(0)>0$  and   the      equality  holds in  a1)   $X(t)\geq 0$ on $[0,2\pi]$
and therefore,   since $f'_t= iX(t) e^{it}$,    $\theta =\arg (f'_t)= t + \pi/2$. Hence if $f$ is homeomorphism,    $\gamma_f$ is convex.
\begin{proof} We  suppose first that   $f$ is harmonic on  $\overline{\mathbb{U}}$ (in general case we can  apply the obtained results on $f_r$, $0<r<1$, and then pass by limit when  $r$ tends $1$).
Set  $z=r e^{it}$. By calculation
$f'_t(z)= i g' re^{it}  +\overline{i h' re^{it}}= i g'(z) z + \overline{i h'(z) z}$. Hence
$f'_t(z)=i \sum_{k=1}^\infty k g_k z^k -i  \overline{\sum_{k=1}^\infty k h_k z^k }$  and  $f'_t(z)=i \sum_{k=1}^\infty k g_k r^ke^{ikt}-
i \overline{\sum_{k=1}^\infty k h_k r^ke^{ikt} }$   and $2 \pi i k  g_k= \int_0^{2\pi} f'_t e^{-ikt}dt$.\\
a) Since   $f'_t= i g' e^{it}  +\overline{i h' e^{it}}$,  $2 \pi i g'(0)= \int_0^{2\pi} f'_t e^{-it}dt$. Hence\\
$2 \pi |g'(0)|\leq  \int_0^{2\pi} |f'_t e^{-it}|dt=L$.\\
Set  $iX(t)=f'_t e^{-it}$,    $X(t)= g' - \overline{h' e^{2i t}}$. Then      (i)  $2 \pi  g'(0)= \int_0^{2\pi} X(t) dt$.\\
If  the equality holds in (i)    and $a= g'(0) >0$, then
$X=X^+$  is a  nonnegative function.
Set   $u=P[X]$  and  $H=h' z^2$. Then   $u= g'- \overline{H}$  and  $u$  is a  nonnegative function.
Hence  ${\rm Im}  ( g')= {\rm Im} (\overline{H})$  and therefore
$g'= -  H + a$, that is  (i6)  $g'= -  h' z^2 + a$.  Since   $X= -  H + a- \overline{H}=a  -{\rm Re}H$, we conclude that $2 H \in\mathcal{P}_a$.\\
Set $M_0(w)= \frac{w}{1+w} $  and $\omega=z^2\nu$.  Then   $2 M_0(w)\in \mathcal{P}_1$  iff  $w\in \mathbb{U}$.\\
It is convenient to suppose for a moment that  $a=1$.
Substitute   $h'= \nu g'$ in  (i6),
we find  $g'=1-z^2\nu g'$ and therefore
$$g'=  \frac{1}{1+ z^2\nu}\quad    \mbox{and}\quad   h'= \frac{\nu}{1+ z^2\nu}.$$
Therefore  $H(z)= M_0(z^2\nu)$  and  $\omega\in {\rm Hol}(\mathbb{U},\mathbb{U})$.\\
Using it one can check first  that
the equality holds in (a3)    in the case $a=1$   and $f(0)=0$ iff  $f=f^{\underline{\nu}}$ and in general       iff $f$ is given by  {\rm(i4)}.

b) For $z \in \overline{\mathbb{U}}$ apply a) on   $f\circ \varphi_z$.
\end{proof}
For the convenience of the reader we first consider  harmonic maps of $\mathbb{U}$ into $\mathbb{R}^3$.
Recall  we will use the following hypothesis in the  sequel   {\rm  ($\underline{H}'_3)$}:  Suppose that  $f=(f^1,f^2,f^3): \mathbb{U}\rightarrow \mathbb{R}^3$  harmonic,   $S=f(\mathbb{U})$
and     the generalized length of  $\partial S$ with respect to $f$, $L=L^+(f)=L^+(f,\partial S)$  is finite.\\
In this setting, let  $F_k$ are holomorphic function in $\mathbb{U}$  such that  $f_k=2 {\rm Re} F_k$, $k=1,2,3$. \\
Then  (A2:)  $f_1+if_2= g +\overline{h}$, where $g=F_1+iF_2$ and $h=F_1-iF_2$.\\
{\rm(II)}   By  $\mathbf{y}=(y_1,y_2,y_3)$  we denote coordinates in $\mathbb{R}^3$  and  for  $\mathbf{y}^0\in \mathbb{R}^3$  we denote the translation  $T_{\mathbf{y}^0}$  defined by
$T_{\mathbf{y}^0}(\mathbf{y})= \mathbf{y} - \mathbf{y}^0$.
We use the following procedure:\\
(I-1)   Set  $p_3(\mathbf{y})=p_3(y_1,y_2,y_3)= (y_1,y_2)$, where by  $\mathbf{y}=(y_1,y_2,y_3)$  we denote coordinates in $\mathbb{R}^3$.
Under  the hypothesis {\rm($\underline{H}_3$)},
it is convenient  to introduce for given  $z_0\in \mathbb{U}$   the tangent plane $Z=Z_{\mathbf{y}^0}$, $\mathbf{y}^0= f(z_0)$.  After rotation we can suppose that \\
($h_2$): $Z$ is  $y_1y_2$-plane  which we can identify  with  $\mathbb{C}$-plane. More precisely there is a   rotation   $R_{\mathbf{y}^0}$ around  $\mathbf{y}^0$  such that    $R=R_Z= T_{\mathbf{y}^0}\circ R_{\mathbf{y}^0}$     maps $Z$ onto  $\Pi=\{ (y_1,y_2,0): y_1,y_2\in\mathbb{R}\}$, with $R_Z(z_0)=0$.
Set  $f^*= R\circ f$,  and   $\gamma^*= R\circ \gamma$.
Then
 $f=f_Z:= p_3\circ R\circ f$  is a  harmonic  function from $\mathbb{U}$ into $\mathbb{C}$.

 Using similar   approach as in  the proof  Theorem \ref{length0} {\rm(ii)}, one can prove:
\begin{prop}
Under  the hypothesis   {\rm ($\underline{H}_3$)},  $2{\rm d=dia}(G)\leq L$.
\end{prop}
\begin{thm}\label{length2}  Suppose that  $f=(f^1,f^2,f^3)$   satisfies  {\rm  ($\underline{H}_3$)}.
Then\\
{\rm a)} {\rm (i)} $ \pi \Lambda_f(0) \leq 2 d$, where  ${\rm d=dia}(S)$.\\
{\rm  b)} $ \pi (1-|z|^2)\Lambda_f(z)\leq 2d$, $z\in \overline{\mathbb{U}}$.

\end{thm}
\bpf
Apply  Theorem \ref{length0}  on  $f_Z$.
\epf
For a fixed $z$,  set  $f_H^z= f_H\circ  \varphi_z-f_H(z)$.
\begin{thm}\label{length3}
Under  the hypothesis {\rm  ($\underline{H}'_3$)},\\
{\rm  b)} If $f$ is  conformal at $z$, then  {\rm (ii)}  $2 \pi (1-|z|^2)|f_x'(z)|\leq L$, $z\in \overline{\mathbb{U}}$.\\
{\rm  b1)}The equality holds in {\rm  b)} for some  $z\in \mathbb{U}$ iff  {\rm(iii)}:  $f(\mathbb{U})$ is  in  the tangent plane $Z=Z_{f(z)}$  and  $f_Z=c f_H^z$   or  $f_Z=c\overline{f_H^z}$, $c\in \mathbb{C}$,
and $H$  satisfy  {\rm(h0)}  with $2 H \in \mathcal{P}_1$  and   $H(z)=0$, where $|c|=|f_x'(z)|(1- |z|^2)$ .
\end{thm}
For a fixed $z$  set $Z=Z_{f(z)}$.  If $f$ is  conformal at $z$,  it is easy to check that   $L=L(f)=L(f_Z)$   and   {\rm  (b2)}:  $2 \pi |(f_Z)_x'(0)|\leq L$.
The equality holds in {\rm  (b2)} iff   $f$ satisfies     {\rm(iii)}.

To get filling about Theorem  \ref{length3}, we give some comments in the  following remark.

\begin{remark}
(i)  Using a similar procedure  one can show  that the corresponding version of  Theorem \ref{length2} and  Theorem \ref{length3}  hold  under hypothesis {\rm ($\underline{H}_3$)}, that is  $m\geq 2$.\\
(ii)  The equality case  in {\rm  b)}.\\
Note that  the equality holds in {\rm  b)} for some  $z\in \mathbb{U}$,  if,  for example, $f(\mathbb{U})$ is in a plane say  $Z$   and  $f=f(z) + R$,
where $R:\mathbb{U}\rightarrow Z$
is a composition of  a rotation in $Z$ around $f(z)$ and homotety wrt $f(z)$.  It is interspersed that the family of extremal maps  is much larger then the family  described
in the previous sentence.
Suppose  that the equality holds in {\rm  b)} for some  $z\in \mathbb{U}$. After rotation we can suppose that $Z=Z_{f(z)}$ is  $y_1y_2$-plane  which we can identify  with  $\mathbb{C}$-plane.  Then  $F_1'(z)= iF_1'(z)$ or  $F_1'(z)= -iF_1'(z)$.
In the case   $F_1'(z)= iF_1'(z)$,   the equality holds in {\rm  b)}
iff $f(\mathbb{U})$ is  in  the tangent plane $Z=Z_{f(z)}$  and   {\rm(i5)}:  $f=c f^{\underline{\nu}} \circ  \varphi_z +c_1$,
where $c,c_1\in \mathbb{C}$,   and   $\nu$  satisfy  {\rm(i3)}  with  $\nu(z)=0$.
In the case   $F_1'(z)= -iF_1'(z)$ we leave the reader to state the corresponding statement.
\end{remark}
\begin{proof}
In particular  if {\rm  ($\underline{H}_3$)}(for dimension $m=3$) holds \footnote {{\rm  ($\underline{H}_3$)}:  $f$  is  continuous on  $\overline{\mathbb{U}}$,    harmonic  on $\mathbb{U}$  and  $\gamma_f$ is a rectifiable curve.}, then the theorem holds.\\
We will prove  the theorem under this assumption. By  application   this case to  $f_r$, $0< r<1$, and letting  $r$ tends to $1$,  one can  get general result.

a) Let  $S=f(\overline{\mathbb{U}})$, and      $M_0=f(0)$.  Since   $f$  is  conformal at $0$,  then \\
(c1):  $f_x'(0)=0$  or  (c2):  $f_x'(0)\times f_y'(0)\neq 0$.\\
In the case (c1),  (i) is clear.
In the case (c2)   there is  the      tangent plane $Z$  of  $S$ at $M_0$.\footnote {Note that   after rotation we can suppose that $Z$ is  $y_1y_2$-plane  which we can identify  with  $\mathbb{C}$-plane  and $f(0)=0$. }
Set $\tilde{f}= p_3\circ f^*$, $\tilde{f}=(f^1,f^2)$
and   $\tilde{\gamma}(t)= \tilde{f}(e^{it})$. Then $\tilde{f}= \tilde{g}+\overline{\tilde{h}}$.

Recall  by notation  in {\rm(II)},
$\Pi$    is   tangent plane  of $S^*=f^*(\overline{\mathbb{U}})$   at $0$,   so that
$$(f^*)^3(x)= o(\tilde{f}(x))=\epsilon(x) \tilde{f}_x'(0)x$$
and therefore  $((f^*)^3)'_x(0)=0$. Hence  $\tilde{f}'_x(0)=f'_x(0)$.
Since   $f$  is  conformal at $0$, then  $\tilde{f}$  is  conformal at $0$,    we can suppose wlg that  $\tilde{h}'_x(0)=0$. Thus, since $R_Z$ is an euclidean isometry  in this case we have \\
(ii1)  $|(\tilde{g})'(0)|= |\tilde{f}'_x(0)|=|f'_x(0)|$.\\
If  $\tilde{L}= |\tilde{\gamma}|$ and  $L^*=|\gamma^*|$ are  lengths of   $\tilde{\gamma}$ and $\gamma^*$  respectively, then by Theorem \ref{length1}a) \\
(ii2)    $2 \pi |(\tilde{g})'(0)|\leq \tilde{L}$.\\
Since  $\tilde{\gamma}$   is the projection of  $\gamma^*$, we first conclude that  (ii3)    $\tilde{L}\leq L^*=L$,
and now  by (ii1),(ii2)  and (ii3), we  get \\
(ii4)   $2 \pi |f'_x(0)|= 2 \pi |(\tilde{g})'(0)|\leq \tilde{L}\leq L^*=L$,\\
which yields  the part a).\\
b)  Apply a) on   $f\circ \varphi_z$.  Note  that   $\tilde{L}\leq L$ with equality iff  $f(\mathbb{U})$ is in a plane.

If   for some  $z\in \mathbb{U}$ the equality holds in b),  then  $\tilde{L}=L$ and the equality holds in b) for  the function $f\circ \varphi_z$ at $0$, that is the equality holds in (b2).
In particular $\tilde{L}=L^*$.  Therefore  $tr(\gamma^*)$  is in a plane $\Pi^*$ parallel to   $\Pi$. By an application of the mean value theorem to  $f^*\circ \varphi_{z_0}$, we conclude that
$0=f^*\circ \varphi_{z_0}(0)$  belongs $\Pi^*$  and therefore  $\Pi^*=\Pi$. Then  $f^*$ is  a  planar mapping and we can apply Theorem \ref{length1}.
\end{proof}

\section{Harmonic  and analytic  disks}\label{ss_gen}
\subsection{Harmonic disks}\label{ss_gen1}
If $f: \mathbb{U} \rightarrow \mathbb{R}^m$  is a  vector  harmonic  on $\mathbb{U}$, we call  $S=f(\mathbb{U})$  a harmonic disk with center at  $f(0)$ (defined by $f$).

\begin{thm}\label{t_incr}
If $f: \mathbb{U} \rightarrow \mathbb{R}^m$  is a  vector  harmonic  on $\mathbb{U}$, then  \\
{\rm (a)} $|f|$  and  $|f'_t|$  are subharmonic.\\
{\rm (b)}   $L(r)$   and   $d(r)$ are increasing in $r\in [0,1)$.
\end{thm}
\bpf
For $z_0\in \mathbb{U}$  and $r>0$ small enough,  by the mean value theorem,  $f(z_0)=  \frac{1}{2\pi}\int_0^{2\pi} f(z_0 + r e^{it}) dt$, and therefore
$$|f(z_0)|=  \frac{1}{2\pi}\left|\int_0^{2\pi}  f(z_0 + r e^{it}) dt\right| \leq   \frac{1}{2\pi}\int_0^{2\pi}  |f(z_0 + r e^{it})|dt .$$
\epf
The following  example  shows how  the boundary behavior
of harmonic mappings may differ from that of conformal mappings. \footnote{According
to the Caratheodory extension theorem, a conformal mapping between
two Jordan domains always extends to a homeomorphism of the closures.}

\bex [\cite{dur}] Let
$l(z)=  \frac{z}{1-z}$,   $s(z)= \frac{1}{2}\ln\frac{1+z}{1-z}$  and  $f(z)=  {\rm  Re} l(z) + i {\rm Im} s(z)$.
Observe that   $f(e^{it})= w_0$   on  $0<t<\pi$ and   $f(e^{it})= \overline{w_0}$   on  $\pi<t<2 \pi$, where $w_0=- \frac{1}{2}+i \frac{\pi}{4}$.
In particular, $f$ collapses the upper and lower semicircles to single points. In
fact, it can be proved that   $l$,  $s$,  and  $f$  map the disk  onto  $S_1= \{ {\rm  Re} w>- \frac{1}{2} \}$, $S_2= \{|{\rm  Im} w| <  \frac{\pi}{4}\}$  and   $S_3= \{ S_1\cap S_2\}$ respectively.
\eex

If the map has no continuous extension to $\overline{\mathbb{U}}$  (in particular   the boundary map  collapses) at first sight the following definitions  seems convenient.
\bdefi
Let  $S\subset \mathbb{R}^m$ be a harmonic disk  defined by $f$.

Let  $h_L^1$   denote the family  of  vector  harmonic function  $f: \mathbb{U} \rightarrow \mathbb{R}^m$  for which  $L^+(f)= \sup \{L(f,r): r\in [0,1) \}<\infty $.
If $L^+(f)<\infty $, then there is a boundary function  $f^*$, but in general  $L^+(f)> L(f^*)$.
\edefi
It seems that  the above described result  hold under each of  the following  hypothesis:

{\rm($\underline{H}'_m$)}: Suppose that  $u=(u^1,u^2,...,u^m): \mathbb{U}\rightarrow \mathbb{R}^m$ is  harmonic,   $S=u(\mathbb{U})$
and     the generalized length of  $\partial S$ with respect to $u$, $L=L^+(u)=L^+(u,\partial S)$  is finite.
If {\rm($\underline{H}_m$)} holds then   the generalized length  $L$  is reduced  to the  length of  $\partial S$.\\

We plan in a forthcoming paper to consider the above discussed results  in connections with  hypothesis {\rm($\underline{H}'_m$)}.

Here we only show  that  that the corresponding version of  Theorem  \ref{length0} holds under hypothesis {\rm(H0)}.
\begin{thm}\label{length+0}
Let $f= g+\overline{h}$ be  complex valued harmonic   in   $\mathbb{U}$  which satisfies  {\rm($\underline{H}_m$)}. \\
Then \\
{\rm(i)} $\pi \Lambda_f(0)\leq 2d$.\\
{\rm(ii)} $2d \leq  L^+(f)$.\\
For function
$$u_d(z)= \frac{d}{\pi }\arg \frac{1+z}{1-z}$$
the  equality holds in {\rm(i)}.
\end{thm}
It is interesting that  $p(z)=d \cdot x/2$  is not extremal for the inequality   (i).
\bpf
Set   $\hat{e}_1= df_0(e_1)$. We can suppose that  $f(0)=0$   and using rotations   that   $\Lambda_f(0) =|\hat{e}_1|$ and   $\hat{e}_1= df_0(e_1)=k e_1$, where  $k=\Lambda_f(0)$.
Set  $p(w)= u$  and $F= p \circ f$.  Then  $p(G)$  is an interval of length equal or  less then $d$,
and  by Proposition  \ref{propMM0},   $\pi |\nabla F(0)|\leq 2d$. Since  $|\nabla F(0)|=\Lambda_f(0)$,  we get the first inequality  of  (i).
We   leave to the reader    to show that  $2d \leq  L$.\\
 For  $0<r<1$  set   $G_r=f_r(\mathbb{U})$,  $|\gamma_{f_r}|=L_f(r)$  and  denote with  $d_r$ the  diameter of $G_r$.\\
By   Theorem  \ref{length0},   (iii)  $\pi r \Lambda_f(0)\leq 2d(r)\leq L(r)$,   $0<r<1$.
Since  $|f'_t|$ and $|f|$ are subharmonic  $L_r$ and  $d_r$ are increasing functions in $r \in[0,1)$. Hence by letting $r$ to $1$  in (iii),  one can prove the result.
\epf

\subsection{Harmonic  and analytic  disks}
For $\mathbf{z}=(z_1,z_2,...,z_m)\in \mathbb{C}^m$, set    ${\rm Re}\,\mathbf{z}= ({\rm Re} z_1,{\rm Re} z_2,...,{\rm Re} z_m) $  and  ${\rm Im}\,\mathbf{z}= ({\rm Im} z_1,{\rm Im} z_2,...,{\rm Im} z_m)$.
Recall  we will use the following hypothesis in the  sequel   {\rm($H_m$)}:  Suppose that  $u=(u^1,u^2,...,u^m): \mathbb{U}\rightarrow \mathbb{R}^m$ is  harmonic,   $S=u(\mathbb{U})$
and     the generalized length of  $\partial S$ wrt $u$, $L=L(u)=L_-(u,\partial S)$  is finite.\\
(B0)In this setting, there are    are holomorphic functions  $F_k$  in $\mathbb{U}$  such that  $u_k= {\rm Re} F_k$, $k=1,2,3,...,m$. Set  $F=(F^1,F^2,...,F^m)$.
We say shortly that  holomorphic function  $F$ is associated  to $u$.
Then $u'_x=\frac{1}{2}(F'_x + \overline{F'_x})= {\rm Re} F'(z)$  and  therefore\\
(B)     $u'_x- i u'_y= F'$.\\
If $f:G \rightarrow \mathbb{R}^2$,  recall  then \\
(B1:)  $f_1+if_2= g +\overline{h}$, where $g=(F_1+iF_2)/2$ and $h=(F_1-iF_2)/2$.\\
(B2:) If $p=f_z$  and  $q=f_{\bar{z}}$, then  $p=f_z=g_z=g'$,  $q=f_{\bar{z}}=\bar{h'}$,  $J_f= {\rm Re}(i\overline{F_1'}F_2')$ and \\
$4 |p|^2= |F'|^2 + 2 J_f$, $4 |q|^2= |F'|^2 - 2 J_f$,   $2(|g'|^2 + |h'|^2)=  |F'|^2$   and\\
if  $|g'|\geq |h'|$, then  $2 |g'|\geq |F'|\geq 2|h'|$.\\
(B3:)  If in addition $u$ is conformal at $0$, then  $h'(0)=0$, and  $|F'(0)|= \sqrt{2} |g'|$.
\begin{thm}\label{tH0}
Suppose  the hypothesis   {\rm($\underline{H}_m$)}.   \footnote{ {\rm($\underline{H}_m$):}   $u: \mathbb{U} \rightarrow \mathbb{R}^m$    is   continuous on  $\overline{\mathbb{U}}$  and   harmonic  on $\mathbb{U}$, and   $\gamma=\gamma_u$  is rectifiable.}\\
{\rm d1)} Then      there  is  a  holomorphic function  $F:\mathbb{U} \rightarrow \mathbb{C}^m$
such that  $u={\rm Re}F $;     and in this setting   $\pi |F'(0)|\leq L$,   where  $L= |\gamma|$ is length of   $\gamma$.\\
{\rm d2)} Then  {\rm (i):}    $2 \pi |D_zu(0)|\leq L$.\\
{\rm d3)}  If in addition $u$ is conformal at $0$, then   $2\pi \Lambda_u(0)\leq L$.
\end{thm}
\bpf
Since   $2 u'_t= F'(z) i e^{it} +  \overline{i z F'(z)}$,  $2 \pi\hat{F}(0)= \int_0^{2\pi} F'(z)dt$,  $F'(0)= \hat{F}(1)$  $\pi i F'(0)= \int_0^{2\pi} u'_t e^{-it}dt$, we find   $\pi |F'(0)|\leq \int_0^{2\pi} |u'_t e^{-it}|dt=L$.
Then  $\pi |F'(0)|\leq L$, and since $2D_zu(0)=F'(0)$, we get (i).\\
Using similar   approach as in  the proof of  Theorem \ref{length2}(the  procedure  described in (I-1)) and apply Theorem \ref{length1}(planar case), one can prove  d3).
\epf
\begin{remark}
By (B),  since  $u'_x,u'_y \in \mathbb{R}^m$,    $|F'|^2=|u'_x- i u'_y|^2= |u'_x|^2+ |u'_y|^2$.
Since  $2u'_t= i z F'(z)  +  \overline{i z F'(z)}$,   $L\leq  |F'|_1$.

Question 1. What is relation between $L$ and   $|F'|_1$?

Note that   $2D_zu(0)=F'(0)$  in $\mathbb{C}^m$. If $m=2$  then $D_zu= g'$  in $\mathbb{C}$. Here we need to be careful  because we identify $(u_1,u_2)\in \mathbb{C}^2$ with  $u_1+i u_2\in \mathbb{C}$(but the corresponding norms  in $\mathbb{C}^2$ and  $ \mathbb{C}$  are not equal in general). Therefore  $2 |g'| \neq  |F'|$ in general.
(B3)  shows that the estimat (i) in  {\rm d2)} is not optimal in general.

Question 2. Can we modify our procedure to get an optimal estimate?
\end{remark}
\section{area estimate}\label{se_area}
We advise the reader to recall Definition \ref{def2}.

\begin{thm} Suppose  that $u$   satisfies  the hypothesis   {\rm($\underline{H}_m$)}.
Then \\  {\rm (i1)}:  $4 \pi A(S)\leq L^2$, where   $A=A(u)$  and   $L=L(u)$. \\
{\rm(I1)}  If in addition to  {\rm($H_m$)}  we suppose that \\
{\rm(h2):}  $A(S)$  is finite  and \\
{\rm(h3):}  $u$ is  $\underline{K}$-quasiconformal  and  $F=\underline{F}$ is   a corresponding  holomorphic function  associated  to $u$, then
$$D[u]= \underline{D}[F]=\pi \bigl(\sum_{k=1}^\infty (k |\hat{F}(k)|^2) \bigr) \leq 2  \underline{K }A(S).$$

{\rm (I2)}  In particular  under  {\rm ($H_m$)}  and  {\rm (h2)},\\
{\rm (i2):}  $|D_1 u(0)|^2 + |D_2 u(0)|^2 \leq 2 \underline{K }A(S)$ with equality iff \\
{\rm(i3):}   $u(z)= ax +by$,  where  $a=D_1 u(0)$,  $b=D_2 u(0)$   with $ \underline{K }= \frac{|a|^2 + |b|^2 }{J}$, where $J=\sqrt{|a|^2|b|^2  -  (a\cdot b)^2 }$. In the case  {\rm (i2)},  $u(\mathbb{U})$   is     a planar domain bounded by an ellipse.
\end{thm}
\bpf
Since the Gaussian curvature of $S$ is negative, by a version of  isoperimetric inequality
(see for example  Theorem 3.4 \cite{KalMarMat}),   we  get (i1).
Set  $J_u=\sqrt{EG -F^2}$.  Then,  by {\rm(h3):}
$|F'|^2\leq \underline{K }J_u$  on  $\mathbb{U}$.
By (A) we have  $D[u]= \underline{D}[F]$  and hence by Parseval's formula we get   {\rm (I1)}.
If equality holds in   {\rm i2)}   then  $\hat{F}(k)=0$ for $k>1$   and therefore    $F(z)= c z$, where  $c=a+ib$,  $a,b\in \mathbb{R}^m$. Hence  we get (i3).

\epf

\begin{thm}
{\rm d1)}   If in addition to  {\rm ($H_m$)}  we suppose that {\rm(h2):}  $A(S)$  is finite  and   {\rm(h4):}  $u$  is  conformal, then
$$A=\int_{\mathbb{U}} |D_1u|^2 dx dy = \frac{\underline{D}[F]}{2}=\frac{\pi}{2} \bigl(\sum_{k=1}^\infty (k |\hat{F}(k)|^2) \bigr).$$
{\rm d2)}  In particular,    $\pi \Lambda^2_u(0)\leq D[u]$  with equality iff  {\rm (ii):} $\gamma=\gamma_u$  is a circle  given by  $u_k=a_kx -b_ky$, $k=1,2,3,...,m$,   where   $|a|= |b|$  and $a\cdot b=0$.\\
{\rm d3)}   $2\pi \Lambda_u(0)\leq L$   with equality iff (ii).
\end{thm}
Examples   $u(z)= z + \overline{z}$ and  $u_n(z)=n z + \overline{z}/n$  show that  {\rm i3)} is not true in general without  hypothesis that the mapping is qc.
\bpf
By (A) we have  $2A=D[u]= \underline{D}[F]$  and hence by Parseval's formula we get   {\rm d1)}.
If equality holds in   {\rm d2)}   then  $\hat{F}(k)=0$ for $k>1$   and therefore    $F(z)= c z$, where  $c=a+ib$,  $a,b\in \mathbb{R}^m$,
and since  $u$  is  conformal   at $0$  {\rm (ii)}  holds.\\
By the isoperimetric inequality  $2 \pi   D[u]=4 \pi A \leq L^2$  and therefore  {\rm d2)}  implies  {\rm d3)}.
\epf

\section{Appendix}


Let $f:\Omega \rightarrow f(\Omega)$ be a $C^{1}$-diffeomorphism. We write  $df=pdz +q d\overline{z}$, where  $p=f_{z}$ and $q=f_{\overline{z}}$.\\
$J_{f} = {| f_{z} |}^{2} -{| f_{\overline{z}} |}^{2}$.\\
Let  $f$ be a diffeomorphism in a
neighborhood $U$ of a point $z_0$. Then  $f$  is  orientation
preserving mapping in $U$ if and only if  $J_{f}(z_{0})>0$.\\
If  $f$  is  orientation
preserving mapping in $U$ at  $z_0$, then  $df$  maps  the tangent space  $T_{z_0}$ into  $T_{w_0}$, where  $w_0=f(z_0)$, and  circles $K_r$ with center at $z_0$   of radius $r$  onto
ellipses $E_r$ with center at $w_0$  and  with major axis of length  $\Lambda_f r$  and minor   axis of length  $\lambda_f r$.
The dilatation (or distortion) at $z_{0}$ is defined to be
\begin{equation}
D_{f}:= \frac{| f_{z} | + | f_{\overline{z}} | }{| f_{z}| - |
f_{\overline{z}} |} \geq 1.
\end{equation}

The complex dilatation at $z_{0}$ is

\begin{equation}
\mu _{f} = \frac{f_{\overline{z}}}{f_{z}}.
\end{equation}
It is often more convenient to consider
\begin{equation*}
d_f =|\frac{f_{\overline{z}}}{f_z}|.
\end{equation*}
The dilatation and distortion are related by
\begin{equation*}
D_{f} = \frac{1+ | \mu _{f} |}{1 - | \mu_{f} |}.
\end{equation*}
Let $f\in C^{1}$ be  orientation
preserving mapping. Then $f$ is conformal iff  $q=f_{\overline{z}} \equiv
0$ (Cauchy-Riemann equations).
If $f$ is conformal, $D_{f}=1$ and  $q=0$, so $df=p dz$
maps circles to circles.
\begin{definition}[Grotzsch analytic definition for regular
mappings]\label{def.qc.1}
Let $f:\Omega \rightarrow \mathbb{C}$ be a diffeomorphism. We say
that $f$ is a \emph{quasiconformal map} if $D_{f}(z)$ is bounded
in $\Omega $. We say $f$ is a K-quasiconformal map if
$D_{f}(z)\leq K$ for all $z\in \Omega $.
\end{definition}

$K(f)= ess \sup_{z \in \Omega}  D_{f}(z)$  is called the coefficient of quasi-conformality (or linear dilatation) of $f$  in the domain $\Omega$.
\bdefi
b4) For   a planar  domain $D$  and   $C^1$ mapping   $u: D \rightarrow \mathbb{R}^m$,  set  $S= u(D)$,  $$K_*(f,z)=  \frac{E+G}{2 J_u}$$   and    $K_*(f)= ess \sup_{z \in D} K_*(f,z)$  which  is called the coefficient of quasi-conformality (or linear dilatation) of $f$  in the domain $D$.
\edefi
\begin{remark}\label{rmQc3}
If $S$  is in a plane  and   $K$ the standard coefficient of quasi-conformality,  then $K_*= \frac{K^2 +1}{2 K}$,that is   $K= K_* + \sqrt{K_*^2-1}$,  where   $K_*= K_*(f)$ and  $K= K(f)$.
Motivated by this  we
give an alternative definition: $u$ is $K$-qc if   $E+F \leq \underline{K} J_u$, where  $\underline{K}= K + \frac{1}{K}$; in planar case  this definition  is
reduced to the standard definition of coefficient of quasi-conformality.
\end{remark}

If $D,G $ are domains in $\mathbb{R}^n$,  by  ${\rm Har}(D,G)$   denote the family of all vector   valued  harmonics maps
$f$ from  $D$ into $G$.

\bdefi [${\rm Har}(p)$,${\rm Har}_c(p)$]\label{}
 For $p\in \mathbb{B}$, let ${\rm Har}(p)= {\rm Har}(\mathbb{B} ,\mathbb{B}; p)$ (respectively ${\rm Har}_c(p)$)   denote the family of all vector   valued  harmonics maps
$f$ from $\mathbb{B}$ into itself with $f(0)=p$ (respectively  which are conformal at  $0$  respectively).

Set $L_h(p)= \sup \{|f'(0)|: f\in {\rm Har}(p) \}$ and $K_h(p)=\frac{L(p)}{\sqrt{1-|p|^2}}$,  $L_c(p)= \sup \{|f'(0)|: f\in {\rm Har}_c(p) \}$ and $K_c(p)=\frac{L_c(p)}{1-|p|^2}$.

For planar domains $D$ and $G$  and given    $z\in D$ and $q \in  G$   denote by    $L_h(z,p; D,G)= \sup \{|f'(z)|\}$,  where  the supremum is taken over all       $f \in {\rm  Har}(D,G)$ with $f(z)=p$. If $D=\U$  we write ${\rm  Har}(G)$ instead of   ${\rm  Har}(\U,G)$   and if in addition  $z=0$, we write simply  $L_h(p,G)$ (or   $L_{\rm har}(p,G)$ ) and if in addition $G=\U$,   $L_h(p)$.
\edefi

\begin{problem}[Extremal] \label{prob1}
For given $p\in \mathbb{B}$
find $K_h(p)$  and   $K_c(p)$.
\end{problem}

For given $p\in \mathbb{B}$, find   $\sup \{|f'(p)|: f\in {\rm Har}(\mathbb{B} ,\mathbb{B})\}$.
For given $p,q\in \mathbb{B}$, find   $\sup \{|f'(p)|: f\in {\rm Har}(\mathbb{B} ,\mathbb{B}), f(p)=q\}$.

The editors of JMAA  paid   my attention to   \cite{Khav}  and the book \cite{Maz}.

{\it Acknowledgement}. We are indebted  to Shi Qingtian, who has been  reading  very carefully  several versions,  for useful   discussions  and useful comments which improved the exposition.


\end{document}